\documentclass[11pt,twoside]{amsart}
\newtheorem{thm}{Theorem}[section]
\newtheorem{cor}[thm]{Corollary}
\newtheorem{lem}[thm]{Lemma}

\newtheorem{defn}[thm]{Definition}
\newtheorem{rem}[thm]{\bf{Remark}}

\numberwithin{equation}{section}
\def\pn{\par\noindent}

\begin{document}

%--------------------------------------------------
%%Don not change any thing in this part
\leftline{ \scriptsize \it Bulletin of the Iranian Mathematical
Society  Vol. {\bf\rm XX} No. X {\rm(}200X{\rm)}, pp XX-XX.}

\vspace{1.3 cm}

%----------------------------------------------------------------------------
\title
{On a Conjecture of a Bound for the Exponent
of the Schur Multiplier of a Finite $p$-Group}
\author{B. Mashayekhy*, A. Hokmabadi and F. Mohammadzadeh}

\thanks{{\scriptsize This research was supported by a grant from Ferdowsi University of Mashhad; (No. MP87150MSH).\\
%This research is supported by the research Council of Sharif University of Technology\\
\hskip -0.4 true cm MSC(2000): Primary: 20C25, 20D15; Secondary: 20E10, 20F12.
%$\dag$Her research supported by a grant of the research council of the University of Tehran.\\
Received: April 17, 2009, Accepted: August 1, 2010.\\
$*$Corresponding author
\newline\indent{\scriptsize $\copyright$ 2010 Iranian Mathematical
Society}}}

\maketitle

%-----------------------------------------------------------------------------
%This part will be filled in by BIMS
\begin{center}
Communicated by\;
\end{center}
%----------------------------------------------

\begin{abstract}  Let $G$ be a $p$-group of nilpotency class
$k$ with finite exponent $\exp(G)$ and let $m=\lfloor\log_pk\rfloor$. We show that $\exp(M^{(c)}(G))$ divides
$\exp(G)p^{m(k-1)}$, for all $c\geq1$, where $M^{(c)}(G)$ denotes the c-nilpotent multiplier of $G$. This implies that
$\exp( M(G))$ divides $\exp(G)$ for all finite $p$-groups of class at most $p-1$. Moreover, we show that our result
is an improvement of some previous bounds for the exponent of
$M^{(c)}(G)$ given by M. R. Jones, G. Ellis and P. Moravec in some cases.
\end{abstract}

\vskip 0.2 true cm

%-----------------------------------------------------------------------------

\pagestyle{myheadings}
\markboth{\rightline {\scriptsize  B. Mashayekhy, A. Hokmabadi and F. Mohammadzadeh}}
         {\leftline{\scriptsize On a Conjecture of a Bound for the Exponent of the Schur Multiplier}}

\bigskip
\bigskip

%-----------------------------------------------------------------------------
%-----------------------------------------------------------------------------

\vskip 0.4 true cm

\section{\bf Introduction and Motivation}

\vskip 0.4 true cm
Let a group $G$ be presented as a quotient of a free group $F$ by
a normal subgroup $R$. Then the $c$-nilpotent multiplier of $G$ (the Baer invariant of $G$
with respect to the variety of nilpotent group of class at most $c$) is
defined to be $$M^{(c)}(G) =\frac{R \cap \gamma _{c+1}(F)}{[R,\
_cF]},$$where $[R,\ _cF]$ denotes the commutator subgroup $[R,
\underbrace{F,...,F}_{c-times}]$ and $ c \geq 1$. The case $c=1$
which has been much studied is the Schur multiplier of $G$, denoted by $M(G)$.
When $G$ is finite, $M(G)$ is isomorphic to the second cohomology group $H^{2}(G, {\mathbb C}^{*})$
(see G. Karpilovsky [6] and C. R. Leedham-Green and S. McKay [8] for further
details).

It has been interested to finding a relation between the exponent
of $M^{(c)}(G)$ and the exponent of $G$. Let $G$ be a finite
$p$-group of nilpotency class $k\geq2$ with exponent $\exp(G)$. M. R. Jones [5] proved that
$\exp(M(G))$ divides $\exp(G)^{k-1}$. This has been improved by G. Ellis [3] who showed that
$\exp( M^{(c)}(G))$ divides $\exp(G)^{\lceil
k/2\rceil}$, where $\lceil k/2\rceil$ denotes the smallest integer
$n$ such that $n \geq k/2$. For $c=1$, P. Moravec [11] showed that
$\lceil k/2\rceil$ can be replaced by $2\lfloor \log_2{k}\rfloor$
which is an improvement if $k\geq 11$.

In this paper we will show that if $G$ is a finite exponent $p$-group of class $k\geq 1$, then $\exp(M^{(c)}(G))$ divides $\exp(G)p^{m(k-1)}$, for all $c\geq1$, where $m=\lfloor\log_pk\rfloor$.
Note that this result is an
improvement of the results of Jones, Ellis and Moravec if $\lfloor\log_pk\rfloor(k-1)/k<e$, $\lfloor\log_pk\rfloor(k-1)/\lceil k/2\rceil-1<e$,
$\lfloor\log_pk\rfloor(k-1)/2\lfloor\log_2k\rfloor-1<e$, respectively, where $\exp(G)=p^e$.

It was a longstanding open problem as to wether $\exp(M(G))$ divides $\exp(G)$ for every finite group $G$. In fact it was conjectured that the exponent of the Schur
multiplier of a finite $p$-group is a divisor of the exponent of
the group itself. I. D. Macdonld and J. W. Wamsley [1]
constructed an example of a group of order $2^{21}$ which has exponent 4, whereas its Schur
multiplier has exponent 8, hence the conjecture is not true in
general. Also Moravec [12] gave an example of a group of order $2048$ and nilpotency class 6 which has exponent 4 and multiplier of exponent 8.
He also proved that if $G$ is a group of exponent 4, then
$exp(M(G))$ divides 8. Nevertheless, Jones [5] has shown that the
conjecture is true for $p$-groups of class 2 and emphasized that
it is true for some $p$-groups of class 3. S. Kayvanfar and M. A.
Sanati [7] have proved the conjecture for $p$-groups of class 4
and 5, with some conditions. A. Lubotzky and A. Mann [9] showed that the conjecture is true for powerful $p$-groups.
The first and the third authors [10] showed that the conjecture is true for nilpotent multipliers of powerful $p$-groups.
Finally, Moravec [11, 12] showed that the conjecture is true for metabelian groups of exponent $p$, $p$-groups with potent filtration and $p$-groups of maximal class. Note that a consequence of our result shows that the conjecture is true for all
finite $p$-groups of class at most $p-1$.

\vskip 0.4 true cm

\section{\bf {\bf \em{\bf Preliminaries}}}

\vskip 0.4 true cm

In this section, we are going to recall some notions we will use
in the next section.
\vskip 0.4 true cm
\begin{defn} (M. Hall [4]).
Let $X$ be an independent subset of a free group, and select an
arbitrary total order for $X$. We define the basic commutators on
$X$, their weight \textit{wt}, and the ordering among them as
follows:

(1) \ The elements of $X$ are basic commutators of weight one,
ordered according to the total order previously chosen.

(2) \ Having defined the basic commutators of weight less than $n$,
the basic commutators of weight $n$ are the $c_k=[c_i,c_j]$, where:

(a) \ $c_i$ and $c_j$ are basic commutators and $wt(c_i)+wt(c_j)=n$,
and

(b) \ $c_i>c_j$, and if $c_i=[c_s,c_t]$, then $c_j\geq c_t$.

(3) \ The basic commutators of weight $n$ follow those of weight
less than $n$. The basic commutators of weight $n$ are ordered among
themselves lexicographically; that is, if $[b_1,a_1]$ and
$[b_2,a_2]$ are basic commutators of weight $n$, then $[b_1,a_1]\leq
[b_2,a_2]$ if and only if $b_1<b_2$ or $b_1=b_2$ and $a_1<a_2$.
\end{defn}
\vskip 0.4 true cm
\begin{lem} (R. R. Struik [13]). Let $x_1, x_2,
... , \ x_r$ be any elements of a group and let
$\upsilon_1,\upsilon_2,...$  be the sequence of basic commutators
of weight at least two in the $x_i$'s, in ascending order. Then
$$(x_1 x_2 ... x_r)^{\alpha}=x_{i_1}^{\alpha}x_{i_2}^{\alpha} ...x_{i_r}^{\alpha} \upsilon_1^
{f_1(\alpha)}\upsilon_2^{f_2(\alpha)}...\upsilon_i^{f_i(\alpha)}...
\ ,$$ where $\{i_1, \ i_2, ...,\ i_r \} = \{1, 2, ..., r \}$, $\alpha$ is a nonnegative integer and
$$ f_i(\alpha)=a_1{\alpha \choose 1}+ a_2{\alpha \choose 2}+...+a_{w_i}{\alpha \choose
w_i}, \ \ \ \ ( I ) \ $$ with $a_1, ..., a_{wi} \in \mathbf{Z} $ and $w_i$ is the
weight of $\upsilon_i$ in the $x_i$'s.
\end{lem}
\vskip 0.4 true cm
\begin{lem} (Struik [13]). Let $\alpha$
be a fixed integer and $G$ be a nilpotent group of class at most
$k$. If $b_1, ..., b_r \in G$ and $r<k$, then
$$[b_1,...,b_{i-1},b_i^{\alpha},b_{i+1},...,b_r]=[b_1,...,b_r]^{\alpha}
\upsilon_1^{f_1(\alpha)}\upsilon_2^{f_2(\alpha)}...,$$ where
$\upsilon_i$'s are commutators in $b_1,...,b_r$ of weight strictly
greater than $r$, and every $b_j$, $1\leq j \leq r$, appears in
each commutator $\upsilon_i$, the $\upsilon_i$'s listed in
ascending order. The $f_i(\alpha)$'s are of the form (I), with $a_1, ..., a_{w_i} \in
\mathbf{Z}$ and $w_i$ is the weight of $\upsilon_i$ ( in the $b_j$'s )
minus $(r-1)$.
\end{lem}
\vskip 0.4 true cm
\begin{rem} Outer
commutators on the letters $x_{1}, x_{2},\ldots,x_{n},\ldots$ are
defined inductively  as follows:

 The letter $x_i$ is an outer commutator
word of weight one. If $u=u(x_1,\ldots,x_s)$ and
$v=v(x_{s+1},\ldots,x_{s+t})$ are outer commutator words of weights
$s$ and $t$, then
$w(x_1,\ldots,x_{s+t})=[u(x_1,\ldots,x_s),v(x_{s+1},\ldots,x_{s+t})]$
is an
outer commutator word of weight $s+t$ and will be written $w=[u,v]$.

It is noted by Struik [13]
that Lemma 2.3 can be proved by a similar method if
$[b_1,..,b_{i-1},b_i^{\alpha},b_{i+1},...,b_r]$ and
$[b_1,...,b_r]$ are replaced with outer commutators.
\end{rem}

By a routine calculation we have the following useful fact.
\vskip 0.4 true cm
\begin{lem} Let $p$ be a prime number and $k$ be
a nonnegative integer. If $m=\lfloor\log_pk\rfloor$,
then $p^t$ divides $p^{m+t} \choose k$, for all integers $t\geq 1$.
\end{lem}

\vskip 0.4 true cm

\section{\bf {\bf \em{\bf Main Results}}}

\vskip 0.4 true cm

In order to prove the main result we need the following lemma.
\vskip 0.4 true cm
\begin{lem} Let $G$ be a $p$-group of
class $k$ and exponent $p^e$ with a free presentation $F/R$. Then for any $c\geq 1$,
every outer commutator of weight $w > c$ in $F /[R,\ _c F]$ has
an order dividing $p^{e+m(c+k-w)}$, where $m=\lfloor\log_pk\rfloor$.
\end{lem}
\vskip 0.4 true cm
\pn{\bf Proof.} Since $\gamma_{k+1}(F) \subseteq
R$, we have $\gamma_{c+k+1}(F) \subseteq [R,\ _cF]$. Also, for all
$x$ in $F$ and $t \geq 0$ we have $x^{p^{e+t}} \in R$ and hence
every outer commutator of weight $w > c$ in $F$, in which
$x^{p^{e+t}}$ appears, belongs to $[R,\ _cF]$. Now we use inverse
induction on $w$ to prove the lemma. For the first step, $w=c+k$,
the result follows by the above argument and Lemma 2.3. Now assume
that the result is true for all $l>w$. Put $\alpha=p^{e+m(c+k-w)}$
and let $u=[x_1, \dots , x_w]$ be an outer commutator of weight
$w$. Then by Lemma 2.3 and Remark 2.4 we have
$$ [x_1^{\alpha}, \dots , x_w]=[x_1, \dots , x_w]^{\alpha}
\upsilon_1^{f_1(\alpha)} \upsilon_2^{f_2(\alpha)} \dots,$$ where the $\upsilon_i^{f_i(\alpha)}$
are as in Lemma 2.3. Note that $w < w_i = wt(\upsilon_i) \leq c+k
$ modulo $[R,\ _cF]$ and hence $f_i(\alpha)=a_1{\alpha \choose
1}+ a_2{\alpha \choose 2}+...+a_{w_i}{\alpha \choose k_i}$, where
$k_i=w_i - w + 1 \leq c+k-w+1 \leq k$, for all $i\geq 1 $.
Thus Lemma 2.5 implies that $p^{e+m(c+k-w-1)}$ divides the
$f_i(\alpha)$'s. Now by induction hypothesis
$\upsilon_i^{f_i(\alpha)} \in [R,\ _cF]$, for all $i\geq 1 $. On the other hand, since $x_1^{\alpha} \in R$ and $w>c$,
$[x_1^{\alpha}, \dots, x_w] \in [R,\ _cF]$. Therefore $u^{\alpha}
\in [R,\ _cF]$ and this completes the
proof.
\vskip 0.4 true cm
\begin{thm} Let $G$ be a $p$-group of
class $k$ and exponent $p^e$. Let $G = F/R$ be any free
presentation of $G$. Then the exponent of $\gamma_{c+1}(F)/[R,\ _cF]$ divides
$p^{e+m(k-1)}$, where $m=\lfloor\log_pk\rfloor$, for all $c\geq 1$.
\end{thm}
\vskip 0.4 true cm
\pn{\bf Proof.} It is easy to see that every
element $g$ of $\gamma_{c+1}(F)$ can be expressed as $g=y_1 y_2
\dots y_n$, where $y_i$'s are commutators of weight at least
$c+1$. Put $\alpha = p^{e+m(k-1)} $. Now Lemma 2.2 implies the
identity
$$g^{\alpha} = y_{i_1}^{\alpha} y_{i_2}^{\alpha} \dots y_{i_n}^{\alpha}
\upsilon_1^{f_1(\alpha)} \upsilon_2^{f_2(\alpha)} \dots,$$ where $\{i_1, i_2, \dots,i_n \}=\{1,
2, \dots, n \}$ and $\upsilon_i^{f_i(\alpha)}$'s are as in Lemma
2.2. Then the $\upsilon_i$'s are basic commutators of weight at
least two and at most $k$ in the $y_i$'s modulo $[R,\ _cF]$ (note that
$\gamma_{c+k+1}(F) \subseteq [R,\ _cF]$). Thus Lemma 2.5 yields
that $p^{e+m(k-2)}$ divides the $f_i(\alpha)$'s. Hence
$\upsilon_i^{f_i(\alpha)} \in [R,\ _cF]$, for all $i\geq 1 $ and $y_j^{\alpha} \in [R,\ _cF]$, for all $1\leq j \leq n$, by
Lemma 3.1. Therefore we have $g^{\alpha} \in [R,\ _cF]$ and the
desired result now follows.

\vskip 0.4 true cm
Now, we are in a position to state and prove the main result of
the paper.
\vskip 0.4 true cm
\begin{thm} Let $G$ be a $p$-group of
class $k$ and exponent $p^e$. Then $\exp( M^{(c)}(G))$
divides $\exp(G)p^{m(k-1)}$, where $m=\lfloor\log_pk\rfloor$, for all $c\geq 1$.
\end{thm}
\vskip 0.4 true cm
\pn{\bf Proof.} Let $G = F/R$ be any free
presentation of $G$. Then $ M^{(c)}(G) \leq \gamma_{c+1}(F)/[R,\
_cF]$. Therefore $\exp( M^{(c)}(G))$ divides
$\exp(\gamma_{c+1}(F)/[R,\
_cF])$. Now the result follows by Theorem 2.3.

Note that the above result improves some previous bounds for the exponent of $M(G)$ and $M^{(c)}(G)$ as follows.\\
Let $G$ be a $p$-group of class $k$ and exponent $p^e$, then we have the following improvements.

$(i)$ If $\lfloor\log_pk\rfloor(k-1)/k<e$, then $\exp(G)p^{\lfloor\log_pk\rfloor(k-1)}<\exp(G)^{k-1}$. Hence in this case our result is an improvement of Jones's result [5].
In particular our result improves the Jones's one for every $p$-group of exponent $p^e$ and of class at most $p^e-1$.

$(ii)$ If $\lfloor\log_pk\rfloor(k-1)/\lceil k/2\rceil-1<e$, then $\exp(G)p^{\lfloor\log_pk\rfloor(k-1)}<\exp(G)^{\lceil k/2\rceil}$ which shows that
in this case our result is an improvement of Ellis's result [3]. In particular our result improves the Ellis's one for every $p$-group of exponent $p^e$ and of class $k<p^{e/3}$, for all $k\geq 3$, or of class $k<p^{e/4}$, for all $k\geq 4$.

$(iii)$ If $\lfloor\log_pk\rfloor(k-1)/2\lfloor\log_2k\rfloor-1<e$, then $\exp(G)p^{\lfloor\log_pk\rfloor(k-1)}<\exp(G)^{2\lfloor \log_2{k}\rfloor}$. Thus in this case our result is an improvement of Moravec's result [11]. In particular our result improves the Moravec's one for every $p$-group of exponent $p^e$ and of class $k<e$, for all $k\geq 2$.
\vskip 0.4 true cm
\begin{cor} Let $G$ be a finite
$p$-group of class at most $p-1$, then $\exp( M^{(c)}(G))$ divides
$\exp(G)$, for all $c\geq 1$. In particular $\exp( M(G))$ divides $\exp(G)$.
\end{cor}
\vskip 0.4 true cm
Note that the above corollary shows that the mentioned conjecture on the
exponent of the Schur multiplier of a finite $p$-group holds for all finite $p$-group of class at most $p-1$.
\vskip 0.4 true cm
\begin{rem}  Let $G$ be a finite
nilpotent group of class $k$. Then $G$ is the direct product of
its Sylow subgroups, $G = S_{p_1}\times \dots \times S_{p_n}$.
Clearly $$\exp(G)= \prod^n_{i=1}\exp(S_{p_i}).$$ By a result of
G. Ellis [2, Theorem 5] we have $$ M^{(c)}(G)=M^{(c)}(S_{p_1}) \times \dots
\times M^{(c)}(S_{p_n}).$$ For all $1\leq i \leq n$, put $m_i=\lfloor\log_{p_i}k\rfloor$. Then by Theorem 3.3
we have $$\exp( M^{(c)}(G)) \mid
\exp(G)\prod^n_{i=1}p_i^{m_i(k-1)}.$$
Hence the conjecture on the
exponent of the Schur multiplier holds for all finite nilpotent group $G$ of class at most $Max\{p_1-1, ..., p_n-1\}$, where $p_1, ..., p_n$ are all the distinct prime divisors of the order of $G$.
\end{rem}

%-----------------------------------------------------------------------------
%-----------------------------------------------------------------------------

\vskip 0.4 true cm

\begin{center}{\textbf{Acknowledgments}}
\end{center}
The authors would like to thank the referee for useful comments.\\ \\

%-----------------------------------------------------------------------------
%-----------------------------------------------------------------------------

\bigskip
\bigskip

%{\bf Received: Month xx, 200x}
{\footnotesize \pn{\bf Behrooz Mashayekhy}\; \\ {Center of Excellence in Analysis on
Algebraic Structures},\\ {Department of
Pure Mathematics}, {Ferdowsi University of Mashhad,\\
P. O. Box 1159-91775,} {Mashhad, Iran.}\\
{\tt Email:
mashaf@math.um.ac.ir}\\
{\footnotesize \pn{\bf Azam Hokmabadi}\; \\ {Department of Mathematics}, {Payame Noor University,} {Iran.}\\
{\tt Email:
hokmabadi-ah@yahoo.com}\\
{\footnotesize \pn{\bf Fahimeh Mohammadzadeh}\; \\ {Department of Mathematics}, {Payame Noor University,} {Ahvaz, Iran.}\\
{\tt Email:
fa36407@yahoo.com}\\
\end{document}